\input amstex
\documentstyle{amsppt}
\magnification=\magstep1
\TagsOnRight
\NoBlackBoxes

\hoffset1 true pc
\voffset2 true pc
\hsize36 true pc
\vsize52 true pc

\def\m1{^{-1}}
\def\ov1{\overline}
\def\gp#1{\langle#1\rangle}
\def\cry#1{\operatorname{\frak{Crys}}(#1)}

\catcode`\@=11
\def\subjclass{\let\savedef@\subjclass
 \def\subjclass##1\endsubjclass{\let\subjclass\savedef@
   \toks@{\def\usualspace{{\rm\enspace}}\eightpoint}%
   \toks@@{##1\unskip.}%
   \edef\thesubjclass@{\the\toks@
     \frills@{{\noexpand\rm2000 {\noexpand\it Mathematics Subject
       Classification}:\noexpand\enspace}}%
     \the\toks@@}}%
  \nofrillscheck\subjclass}
\def\paper{\makerefbox\paper\paperbox@}
\def\jour{\makerefbox\jour\jourbox@{\it}}
\def\logo@{}
\catcode`\@=\active

\topmatter

\title
Torsionfree crystallographic groups with indecomposable holonomy
group
\endtitle
\author
V.A.~Bovdi, P.M.~Gudivok, V.P.~Rudko
\endauthor
\dedicatory
\enddedicatory
\leftheadtext\nofrills{ V.A.~Bovdi, P.M.~Gudivok, V.P.~Rudko }
\rightheadtext \nofrills { Torsionfree crystallographic groups}

\abstract Let $K$ be a principal ideal domain, $G$ a finite group, and
$M$ a $KG$-module which as $K$-module is free of finite rank, and on
which $G$ acts faithfully. A {\it generalized crystallographic group}
(introduced by the authors in volume 5 of this Journal) is a group $\frak
C$ which has a normal subgroup isomorphic to
$M$ with quotient $G$, such that conjugation in
$\frak C$ gives the same action of $G$ on $M$ that we started with.
(When $K=\Bbb Z$, these are just the
classical crystallographic groups.)
The $K$-free rank of
$M$ is said to be the {\it dimension} of $\frak C$, the {\it
holonomy group} of $\frak C$ is $G$, and $\frak C$ is called
{\it indecomposable} if $M$ is an indecomposable
$KG$-module.

Let $K$ be either $\Bbb Z$, or its localization $\Bbb Z_{(p)}$ at the
prime $p$, or the ring $\Bbb Z_p$ of $p$-adic integers, and consider
indecomposable torsionfree generalized crystallographic groups
whose holonomy group is noncyclic of order
$p^2$. In Theorem 2, we prove that (for any given
$p$) the dimensions of these groups are not bounded.

For $K=\Bbb Z$, we show in Theorem 3 that there
are infinitely many non-isomorphic indecomposable torsionfree
crystallographic groups with holonomy group the alternating group of
degree 4. In Theorem 1, we look at a cyclic $G$ whose order $|G|$
satisfies the following condition: for all prime divisors $p$ of $|G|$,
$p^2$ also divides $G$, and for at least one $p$, even $p^3$ does. We
prove that then every product of  $|G|$ with a positive integer
coprime to it occurs as the dimension of some indecomposable
torsionfree crystallographic group with holonomy group $G$.

\endabstract

\thanks
The research was supported by OTKA  No.T 037202, No.T 038059 and
No.T034530
\endthanks

\address
\hskip-\parindent {\rm  V.A.~Bovdi\newline
Institute of Mathematics, University of Debrecen\newline
P.O.  Box 12,
H-4010 Debrecen, Hungary\newline
Institute of Mathematics and Informatics, College of Ny\'\i
regyh\'aza\newline
S\'ost\'oi \'ut 31/b, H-4410 Ny\'\i regyh\'aza,
Hungary
\newline
E-mail: vbovdi\@math.klte.hu
\bigskip
\hskip-\parindent V.P.~Rudko, P.M.~Gudivok\newline
Department of Algebra, University of Uzhgorod\newline
88 000 Uzhgorod,  Ukraine
\newline\hskip\parindent
E-mail: math1\@univ.uzhgorod.ua}
\endaddress

\endtopmatter

\document

\head
 1. Introduction
\endhead
Zassenhaus  developed algebraic methods in \cite{11} for studying
the classical crystallographic groups and he pointed out the close
connection between the theory of crystallographic groups and the
theory of integral representations of finite groups. A historical
overview and the present state of the theory of crystallographic
groups   as well as its connections to other branches of
mathematics are described in \cite{9, 10}.

It was shown in  \cite{5, 7} that, in general,  the classification of the
crystallographic groups is a problem of wild  type, in the sense that it is
related to the classical unsolvable problem of describing the canonical
forms of pairs of linear operators acting on finite dimensional vector
spaces. One may therefore focus attention on certain special classes of
crystallographic groups, for example, on groups whose translation group
affords an irreducible (or indecomposable) integral representation of the
holonomy group. In this direction, Hiss and Szczepa\'nski in \cite{6}
proved  that torsionfree crystallographic groups with irreducible
holonomy group do not exist. On the other hand, Kopcha and Rudko in
\cite{7} showed that the problem of describing torsionfree
crystallographic groups with indecomposable cyclic holonomy group of
order $p^n$,  $n\geq 5$, is still of wild type.

The {\it generalized crystallographic groups} which we introduced in
\cite{3} may be defined as follows. Let $K$ be a principal ideal domain,
$G$ a finite group, and $M$ a $KG$-module which as $K$-module is free
of finite rank, and on which $G$ acts faithfully. A generalized
crystallographic group is a group $\frak C$ which has a normal
subgroup isomorphic to $M$ with quotient $G$, such that conjugation in
$\frak C$ gives the same action of $G$ on $M$ that we started with,
and that the extension in question does not split. The $K$-free rank of
$M$ is said to be the {\it dimension} of $\frak C$, and the {\it
holonomy group} of $\frak C$ is $G$. (In the special case of
$K=\Bbb Z$, this matches one of the usual descriptions of
crystallographic groups; for emphasis, we sometimes refer to those as
{\it classical} crystallographic groups.)

In \cite{3}, we looked at indecomposable generalized
crystallographic groups when $K$ is either $\Bbb Z$, or its localization
$\Bbb Z_{(p)}$ at the prime $p$, or the ring $\Bbb Z_p$ of $p$-adic
integers, and either $G$ is a cyclic $p$-group or $p=2$ and $G$ is
non-cyclic of order~4. Retaining this restriction on the choice of $K$ but
allowing $p$ to be arbitrary, we consider here indecomposable
torsionfree generalized crystallographic groups whose holonomy group
is noncyclic of order $p^2$. In Theorem 2, we prove that (for any given
$p$) the dimensions of such groups are not bounded.

For the classical case ($K=\Bbb Z$), we show in Theorem 3 that there
are infinitely many non-isomorphic indecomposable torsionfree
crystallographic groups with holonomy group the alternating group of
degree 4. In Theorem 1, we look at a cyclic $G$ whose order $|G|$
satisfies the following condition: for all prime divisors $p$ of $|G|$,
$p^2$ also divides $G$, and for at least one $p$, even $p^3$ does. We
prove that then every product of  $|G|$ with a positive integer
coprime to it occurs as the dimension of some indecomposable
torsionfree crystallographic group with holonomy group $G$.

\head
2. The main results
\endhead

For the formal statement of our results, we need some terminology.

Let $K$ be a principal ideal domain, $F$ be a field containing $K$ and let
$G$ be a finite group.  Let $M$ be a $K$-free $KG$-module, with a
finite $K$-basis affording a faithful representation $\Gamma$ of $G$ by
matrices over $K$; further, let $FM$ be the $F$-space spanned by this
$K$-basis of $M$, so $M$ becomes a full lattice in $FM$. Let
$\widehat{M}=FM^+/M^+$ be the quotient group of the additive group
$FM^+$ of the linear space $FM$ by the additive group $M^{+}$ of the
module $M$. Then $FM$ is an $FG$-module and $\widehat{M}$ is a
$KG$-module with operations:
$$
g(\alpha m)=\alpha g(m); \qquad g(x+M)=g(x)+M,
$$
where $g\in G$, $\alpha\in F$, $m\in M$, $x\in FM$.

Let $T:G\to\widehat{M}$ be a $1$-cocycle  of $G$ with values in
$\widehat{M}$. Elements of $\widehat M$ being cosets in $FM^+$
modulo
$M^+$, we consider each value $T(g)$ of $T$ a subset of $FM$,
and define the group
$$
\cry{G; M; T}=\{\,(g,x)\mid g\in G, \,x\in T(g)\,\}
$$
with the operation $$(g,x)(g',x')=(gg',\,gx'+x),$$ where $g,g'\in G$,
$x\in T(g)$, $x'\in T(g')$.

The $K$-free rank of $M$ will be called the $K$-dimension of the group
$\cry{G; M; T}$.  When $T$ is not cohomologous to 0, the group
$\cry{G; M; T}$ is called {\it indecomposable} if $M$ is an
indecomposable $KG$-module.

We note  that if $K=\Bbb Z$ and $F=\Bbb R$, then the abstract group
$\cry{G; M; T}$ is isomorphic to a classical crystallographic group.

Let $C^{1}(G, \widehat{M})$ and $B^{1}(G, \widehat{M})$
be the group of the 1-cocycles and group of the 1-coboundaries of $G$ with
values in the module $\widehat{M}$, respectively, so that
$H^{1}(G,\widehat{M}) =C^{1}(G, \widehat{M})/B^{1}(G,\widehat M)$.  The
group $\cry{G; M; T}$ is an extension of $M^{+}$ by  the group $G$.  The
group $\cry{G; M; T}$ is torsionfree  if and only if for each prime order
subgroup $H$ of $G$  the restriction $T|_H$ is not a coboundary.

Using results from \cite{1,2,8} we prove the following two
theorems.

\proclaim {Theorem 1} Let $G$ be a cyclic group of order
$|G|=p_1^{n_1}\cdots p_s^{n_s}$, where $p_1,\ldots ,p_s$ are
pairwise distinct  primes ($n_1\ge 3$ and if $s\geq 2$ then
$n_2\ge 2,\ldots, n_s\ge 2$),  let $m$ be a natural number which
is coprime to $|G|$ and put $d=m\cdot |G|$. Then there
exists a  torsionfree indecomposable classical crystallographic
group  of dimension $d$ with holonomy group isomorphic to
$G$.
\endproclaim

\proclaim {Theorem 2} Let $K$ be either   $\Bbb Z$, $\Bbb
Z_{(p)}$, or $\Bbb Z_{p}$, and let $G\cong C_{p}\times C_{p}$.
Then the $K$-dimensions of the indecomposable torsionfree groups
$\cry{G; M; T}$ are unbounded.
\endproclaim

In \cite{3} we gave a complete description of the indecomposable
torsionfree crystallographic groups whose holonomy group is
$C_2\times C_2$. Moreover, we proved that there  exist at least
$2p-3$ torsionfree crystallographic groups  having cyclic
indecomposable holonomy group of order ${p^2}$. Note that the holonomy
group of an indecomposable torsionfree crystallographic group can
never have prime order. Therefore the following result is also
interesting.

\proclaim {Theorem 3}  There exist infinitely many
non-isomorphic indecomposable torsionfree classical
crystallographic  groups with holonomy group isomorphic to the
alternating group $A_4$ of degree $4$.
\endproclaim

\head
3. Preliminary results and the proof of Theorem 1
\endhead

Let $K=\Bbb Z$, $\Bbb  Z_{(p)}$ or $\Bbb Z_{p}$ as above,
$H_{p^n}=\gp{a\mid a^{p^n}=1}$ be a cyclic group of order ${p^n}$
($n\geq 2$), $\xi_s$ be a $p^s$th primitive root of unity,
$\xi_s^p=\xi_{s-1}$ ($s\geq 1$, $\xi_0=1$).
Define ordered bases $B_i$ for the free $K$-modules
$\frak R_i=K[\xi_i]$ by setting
$$
\align
B_1&=\{\,1,\xi_1,\dots,\xi_1^{p-2}\,\},\\
B_2&=\{\,1,\xi_1,\dots,\xi_1^{p-2},\xi_2,\,\xi_2\xi_1,\dots,
\xi_2^{p-1}\xi_1^{p-2}\,\},
\intertext{and in general (for $i>1$)}
B_i&=B_{i-1}\cup\xi_iB_{i-1}\cup\xi_i^2B_{i-1}\cup\dots\cup
\xi_i^{p-1}B_{i-1},
\endalign
$$
ordered as indicated. Obviously, $|B_i|=\phi(p^i)$ (where $\phi$ is the
Euler function). It is easy to check that each $\frak R_i$ ($i\leq n$) is a
$KH_{p^n}$-module with action defined by
$$
a(\alpha)=\xi_i\cdot\alpha,\quad\quad\quad (\alpha\in \frak R_i).
\tag1
$$
Keep it in mind that $\frak R_i$ is only a $K$-submodule of
$\frak R_{i+1}$, not a $KH_{p^n}$-submodule. We define by
$\widetilde{\xi_i}$ the matrix of multiplication by
$\xi_i$ with respect to the $K$-basis $B_i$ of the ring $\frak R_i$
($i\geq 1$,
$\frak R_0=K$). Note that
$$
\widetilde{\xi_i}^p=E_p\otimes \widetilde{\xi}_{i-1},\quad\quad
\quad (i>1)
$$
where $E_p$ is the identity matrix of degree $p$ and $\otimes$ is
the Kronecker product of matrices.

Let  $\delta_i$ be the matrix representation of $H_{p^n}$ with respect
to the $K$-basis $B_i$ of the $KH_{p^n}$-module $\frak R_i$. From
(1) it follows that
$$
\delta_i(a)=\widetilde{\xi}_i\quad\quad\quad (i\geq 0,\quad
\widetilde{\xi}_0=1)
$$
and $\delta_0,\ldots, \delta_n$ are irreducible
$K$-representations of $H_{p^n}$.

We introduce the following notation. Let $0\leq i\leq j\leq n$.
For each $\alpha\in \frak R_i$ we denote by $\gp{\alpha}_j^i$ the
matrix with $\phi(p^i)$ rows and $\phi(p^j)$ columns in which all
columns are zero except the last which consists  of the
coordinates of $\alpha\in \frak R_i$ in the basis $B_i$.

Now let  $0\leq i<j\leq n$. Thus
$$
\split
 \widetilde{\xi}_i\cdot \gp{\alpha}^i_j&=\gp{{\xi}_i\alpha}^i_j;\\
\gp{\alpha}^i_j&=({\gp{0}^i_{j-1},
\ldots,\gp{0}^i_{j-1}},\gp{\alpha}^i_{j-1} );\\
\gp{\alpha}^i_j\cdot
\widetilde{\xi_j}^k&=(\gp{\alpha_1(k)}^i_{j-1},
\ldots,\gp{\alpha_{p-1}(k)}^i_{j-1},\gp{\alpha_p(k)}^i_{j-1} ),\\
\endsplit
\tag2
$$
where $0\leq k<p$, \quad $0\leq i<j\leq n$,  \quad
$\alpha_{p-k}(k)=\alpha$ and  \quad  $\alpha_s(k)=0$ for
$s\not=p-k$. The matrix $\gp{\alpha}_j^i$ defines an extension of
the $KH_{p^n}$-module $\frak R_i$ by the $KH_{p^n}$-module $\frak
R_j$ in which the following $K$-representation of the group
$H_{p^n}$ is realized:
$$
a\to
\pmatrix
\widetilde{\xi_i}&\gp{\alpha}_j^i\\
0&\widetilde{\xi_j}
\endpmatrix.\tag3
$$

If $\alpha\equiv0 \pmod{p\frak R_i}$, then the $K$-representation
(3) of $H_{p^n}$ is completely reducible and the corresponding
extension of modules is split, i.e.
$$
p\cdot {Ext}_{KH_{p^n}}(\frak R_j, \frak R_i)=0 \quad \quad (i>j).
\tag4
$$

Let $m$  be a natural number and let $A$ be an $m\times m$ matrix
over $K$. Consider the $K$-representations of the cyclic
group $H_{p^n}=\gp{a\mid a^{p^n}=1}$, with $n>2$,  defined by:
$$
\split
\Delta_1=E_m\otimes\delta_0+E_m\otimes\delta_1&:
\quad a\to
\pmatrix
E_m&0\\
0&E_m\otimes\widetilde{\xi_1}
\endpmatrix;\\
&\\
\Delta_2=E_m\otimes\delta_2+\cdots+E_m\otimes\delta_n&:
\quad a\to
\pmatrix
E_m\otimes\widetilde{\xi_2}&&0\\
&\ddots&\\
0&&E_m\otimes\widetilde{\xi_n}
\endpmatrix;\\
&\\
\Gamma_{p,A}^{(m)}=
\pmatrix
\Delta_1&U\\
0&\Delta_2
\endpmatrix&:
\quad a\to
\pmatrix
\Delta_1(a)&U(a)\\
0&\Delta_2(a)
\endpmatrix,\\
\endsplit
$$
where the matrix
$$
U(a)=
\pmatrix
A\otimes\langle1\rangle_2^0&
E_m\otimes\langle1\rangle_3^0&
\ldots&
E_m\otimes\langle1\rangle_n^0\\
E_m \otimes\langle1\rangle_2^1&
E_m\otimes\langle1\rangle_3^1&
\ldots&
E_m\otimes\langle1\rangle_n^1
\endpmatrix
$$
is  called  the {\it intertwining} matrix.

For $n=2$ we define  the following $K$-representation of
$H_{p^2}=\gp{a\mid a^{p^2}=1}$:
$$
\Gamma^{(1)}_{p}:\quad a\to \pmatrix
1&0&\langle1\rangle_2^0\\
&\widetilde{\xi_1}&\langle1\rangle_2^1\\
0&&\widetilde{\xi_2}
\endpmatrix.\tag5
$$

\proclaim {Lemma 1} Let $J_m$ be the lower triangular Jordan block
of degree $m$ with ones in the main diagonal. Then
$\Gamma_{p,J_m}^{(m)}$ (or  $\Gamma^{(1)}_{p}$, respectively) is
an indecomposable $K$-representation of degree $m\cdot|H_{p^n}|$
of the group $H_{p^n}$ ($n\geq 2$)  (or of degree $|H_{p^2}|$ of
the group $H_{p^2}$\quad for $n=2$, respectively).
\endproclaim
\demo{Proof} Representations which
depend  on  matrix parameters in this way were studied in \cite{1,2}.
Using methods and results from these papers,  it is  not
difficult to show for $n>2$ that the
$K$-representations $\Gamma_{p,A}^{(m)}$ and $\Gamma_{p,B}^{(m)}$
are equivalent if and only if
$$
C\m1AC-B\equiv 0\pmod{p},\tag6
$$
for some invertible matrix $C$. Moreover, the $K$-representation
$\Gamma_{p,A}^{(m)}$ is decomposable if and only if there is a
decomposable matrix $B$, which  satisfies (6). In particular, the
$K$-representation $\Gamma_{p,J_m}^{(m)}$ is an indecomposable
$K$-representation of  $H_{p^n}$.

 The case $n=2$ follows from \cite{1}. The lemma
is proved.
\enddemo

\bigskip

Put
$$
\Gamma_{p}^{(m)}=
\cases
\Gamma_{p,J_m}^{(m)} & \text{  for  } n>2, m>1;\\
\Gamma_{p,1}^{(1)} & \text{  for  } n>2, m=1;\\
\Gamma_{p}^{(1)} & \text{  for  } n=2.\\
\endcases
\tag7
$$
\bigskip

\proclaim {Lemma 2} Let $L_p$ be a $KH_{p^n}$-module affording the
$K$-representation $\Gamma_{p}^{(m)}$ of the group $H_{p^n}$
($n\geq 2$) and $v_1,v_2,\ldots,v_t$ be a $K$-basis corresponding
to this representation in $L_p$. Then $Kv_1$ is a
$KH_{p^n}$-submodule in $L_p$, which over $K$ has a direct
complement $L_p'$ with the following $K$-basis
$\{w_i=v_i+\lambda_i v_1 \}$  for some $\lambda_i\in K$ \quad
($i=2,\ldots, t$) and  which is left invariant with respect to the
operator $a^p$.
\endproclaim
\demo{Proof} Let $n>2$. Clearly, $a\cdot v_1=v_1$, i.e. $Kv_1$ is
a  $KH_{p^n}$-submodule in $L_p$. Using (2) it is easy to check
that in the matrix $\Gamma_{p}^{(m)}(a^p)$ the intertwining matrix
$$
U(a^p)=\sum_{t=0}^{p-1}\Delta_1^{p-t-1}(a)\cdot U(a)\cdot \Delta_2^{t}(a)
$$
has the form
$$
U(a^p)=
\pmatrix
J_m\otimes U_{11}&
\ldots&
E_m\otimes U_{1\,n-1}\\
E_m\otimes U_{21}&
\ldots&
E_m\otimes U_{2\,n-1}\\
\endpmatrix,
$$
where
$$
U_{1\,i}=(\gp{1}_i^0,\ldots,\gp{1}_i^0), \quad\quad U_{2\,i}=(
\gp{1}_i^1,\gp{\xi_1}_i^1, \ldots,\gp{\xi_1^{p-1}}_i^1 ) \quad
(i=1,\ldots,n-1).
$$
We  change the basis elements $v_{m+i}$ by
$w_{m+i}=v_{m+i}+v_1$\quad ($i=1,\ldots,p-1$).  Since the sum
$-(v_{m+1}+\cdots+v_{m+p-1})+v_1$ is replaced  by
$-(w_{m+1}+\cdots+w_{m+p-1})+pv_1$, this  changes the first row of
the matrix $U(a^p)$ turning its elements either to zero or to
multiples of  $p$. From (4)\quad (for $i=0$) provides the
possibility of changing the basis elements as
$$
w_{m+i}=v_{m+i}+\lambda_i\cdot v_1\quad\quad  (p\leq i,\quad
\lambda_i \in K),
$$
and so we get a $K$-module $L_p'$, such that $L_p=Kv_1\oplus L_p'$
and $L_p'$ is left invariant with respect to  the operator $a^p$.

For $n=2$ the statement of the lemma is trivial and the lemma is
proved.
\enddemo

\bigskip

For the remainder of  this chapter we suppose that  $K=\Bbb Z$ and
we will consider only classical crystallographic groups. Let
$G=H_{p_1^{n_1}}\times \cdots \times H_{p_s^{n_s}}$ be a
decomposition of the cyclic group of order $|G|=p_1^{n_1}\cdots
p_s^{n_s}$ into the cyclic subgroups $H_{p_i^{n_i}}$ of order
$|H_{p_i^{n_i}}|={p_i}^{n_i}$, where $p_1,\ldots ,p_s$ are
pairwise distinct   primes,  $n_1\ge 3$ and if $s\geq 2$ then
$n_2\ge 2,\ldots,n_s\ge 2$.

Define by $\Gamma^{(m)}$ the   tensor  product of the $\Bbb
Z$-representation $\Gamma_{p_1}^{(m)}$ of the group
$H_{p_1^{n_1}}$ and the $\Bbb Z$-representations
$\Gamma_{p_j}^{(1)}$ of the groups $H_{p_j^{n_j}}$ ($m\in \Bbb
N$,\quad $j=2,\ldots,s$). Then $\Gamma^{(m)}$ is a $\Bbb
Z$-representation of the group $G$ in which
\smallskip
$$
\Gamma^{(m)}(a_1^{t_1},\ldots,a_s^{t_s})=
\Gamma_{p_1}^{(m)}(a_1^{t_1})\otimes
\Gamma_{p_2}^{(1)}(a_2^{t_2})\otimes\cdots\otimes
\Gamma_{p_s}^{(1)}(a_s^{t_s}).
$$
\smallskip

\proclaim {Lemma 3} If $(m,|G|)=1$ then $\Gamma^{(m)}$ is an
indecomposable  $\Bbb Z$-representation of $G$.
\endproclaim
\demo{Proof} Let $\Gamma^{(m)}|_{H_{p_i^{n_i}}}$ be the
restriction of the representation $\Gamma^{(m)}$ to
${H_{p_i^{n_i}}}$. From Lemma 1 it follows that  the degree of
each indecomposable summand of $\Gamma^{(m)}|_{H_{p_i^{n_i}}}$ is
$m\cdot |{H_{p^n_1}}|$ for $i=1$ and $|{H_{p_i^{n_i}}}|$ for
$i>1$, respectively.

If $\Gamma^{(m)}$ is decomposable and $\Gamma$ is a summand in
$\Gamma^{(m)}$, then the degree of $\Gamma$ is  divisible  by the
number $m\cdot |{H_{p_1^{n_1}}}|$ and if $s\geq 2$ then by the
numbers $|{H_{p_2^{n_2}}}|, \ldots, |{H_{p_s^{n_s}}}|$ (see Lemma
1 for the case $K=\Bbb Z_p)$.

Thus from the condition $(m,|G|)=1$ it follows that
$\Gamma=\Gamma^{(m)}$. The lemma is proved.

\bigskip
Now we construct a cocycle of the group $G$. Let $M$ be a $\Bbb
ZG$-module  of the $\Bbb Z$-representation $\Gamma^{(m)}$
affording  the group $G$ and
$$
M=L_{p_1}\otimes_{K}\cdots\otimes_{K}L_{p_s},\tag8
$$
where $L_{p_i}$ is a $\Bbb ZH_{p^n_i}$-submodule for
$\Gamma_{p_i}^{(m)}$ \quad $(i=1,\ldots, s)$. If
$g=a_1^{t_1}\cdots a_s^{t_s}\in G$ and $l=l_1\otimes\cdots\otimes
l_s\in M$, then
$$
g\cdot l=a_1^{t_1}\cdot l_1\otimes\cdots\otimes a_s^{t_s}\cdot
l_s,
$$
where $l_i \in L_{p_i}$, $t_i\in \Bbb Z$\quad  ($i=1,\ldots,s$).

We can suppose that $M\subset \Bbb{R}^d$. Each $\Bbb Z$-basis in
$M$ is also an $\Bbb R$-basis  in $\Bbb R^d$ and  an $\Bbb
R^+/\Bbb Z^+$-basis in $\widehat{M}={\Bbb R^d}^+/M^+$, where
$d=m\cdot |G|=\deg(\Gamma^{(m)})$.
\bigskip

Let $v=v_1^{(1)}\otimes\cdots\otimes v_s^{(1)}$ be the tensor
product of the first $\Bbb Z$-basis elements of the modules
$L_1,\ldots,L_{p_s}$. Obviously,  $a\cdot v=v$. Define a function
$f:G\to \widehat{M}$ by
 $$
f(g)=\big(\frac{t_1}{p_1^{n_1}}+\cdots+\frac{t_s}{p_s^{n_s}}\big)\cdot
v+M,\tag9
$$
where  $g=a_1^{t_1}\ldots a_s^{t_s}\in G$, $t_1,\ldots,t_s\in \Bbb
Z$.

Since  $g_1\cdot v=v$ and $f(g_1\cdot g_2)=f(g_1)+f(g_2)$ for
$g_1, g_2\in G$, we obtain
$$
f(g_1\cdot g_2)=f(g_2)+f(g_1)=g_1\cdot f(g_2)+f(g_1)
$$
and, therefore,  $f$ is a $1$-cocycle of $G$ in $\widehat{M}$. The
lemma is proved.
\enddemo

\proclaim {Lemma 4} The restriction of the cocycle $f$ to any
prime order subgroup  of $G$ is not a coboundary.
\endproclaim
\demo{Proof} Let $1\le i\le s$ and put  $b=a_i^{p_i^{n_i-1}}\in
G$. From Lemma 2 and (8) it follows that the $\Bbb Z$-module $M$
can be decomposed as  $M=\Bbb Zv\oplus M'$, where $\Bbb Zv$ is a
$\Bbb ZG$-module, $M'$ is a $\Bbb Z$-module which is left
invariant with respect to the  operator $a_i^{p_i}$ (moreover with
respect to the operator $b$). Thus $ \widehat{M}=Fv\oplus
\widehat{M'}$ and $b(\widehat{M'})=\widehat{M'}$. If  $z\in
\widehat{M}$, then $ z=\alpha v+z_1$\quad  ($\alpha\in F, z_1\in
\widehat{M'}$). From (9) and from the condition  $b(z_1)\in
\widehat{M'}$, it follows that
$$
f(b)=\frac{1}{p_i}v+M\ne (b-1)z+M
$$
for any  $z\in \widehat{M}$. Therefore, the restriction of the
cocycle $f$ to $\gp{ b }$ is not a coboundary which proves the
lemma.
\enddemo

\demo{Proof of Theorem 1} From Lemma 4 it follows that $\cry{G;
M;T}$ is a torsionfree group. Moreover, according to Lemma 3,
$\Gamma^{(m)}(G)$ is an indecomposable subgroup in $GL(d, K)$,
where $d=m\cdot |G|$ and $(m,|G|)=1$. So the proof is complete.
\enddemo

\head
4. Proof of Theorem 2
\endhead

Let $K=\Bbb Z$, $\Bbb  Z_{(p)}$ or $\Bbb Z_{p}$ as above and let
$\varepsilon=\xi$ be a $p^{th}$ primitive root of unity ($p>2$).
Then $B_1=\{ 1,\varepsilon, \ldots, \varepsilon^{p-2}\}$ is an
$\frak F$-basis in the field $\frak F(\varepsilon)$ and a
$K$-basis  in the ring $K[\varepsilon]$,  where $\frak F$ is the
field of fractions of the ring $K$.

The symbol $\gp{\alpha}$ denotes a $(p-1)$-dimensional column,
which consists of  coordinates of the element $\alpha\in \frak
F(\varepsilon)$ in the basis $B_1$ and $\widetilde{\alpha}$
denotes the matrix of the operator of multiplication by $\alpha$
in the $\frak F$-basis $B_1$ of the field $\frak F(\varepsilon)$.
Clearly, $\widetilde{\varepsilon}\cdot
\gp{\alpha}=\gp{\varepsilon\alpha}$.

The group $G=\gp{a,b}\cong C_p\times C_p$ ( $p>2$) has the
following $p+2$ irreducible $K$-representations, which are
pairwise nonequivalent over the field $\frak F$:
$$
\matrix \gamma_0 :& a & \to &1, & b &\to &1; &\quad & \gamma_1 :&
a &\to & \widetilde{1},& b&\to& \widetilde{\varepsilon};\cr
\gamma_2 :& a &\to & \widetilde{\varepsilon},& b&\to&
\widetilde{1};&\quad & \gamma_3 :& a &\to &
\widetilde{\varepsilon},& b&\to& \widetilde{\varepsilon};\cr
\rho_i :& a &\to &\widetilde{\varepsilon},& b&\to&
\widetilde{\varepsilon}^i,&& \cr
\endmatrix\tag10
$$
where $i=2,\ldots,p-1$ and $\widetilde{1}=E_{p-1}$ is the identity
matrix of degree $p-1$.

\smallskip Put $\tau=\rho_{p-1}\oplus \cdots \oplus \rho_{2}\oplus
\gamma_3\oplus \gamma_2\oplus \gamma_1$.  Define  the following
$K$-representation of the group $G=\gp{a,b}$:
$$
 \Gamma_0: \quad a\to  \left(\smallmatrix
\tau(a)& U(a)\\
0                                                                 & \gamma_0(a)\\
\endsmallmatrix\right),
\quad\quad\quad\quad\quad\quad
 b\to  \left(\smallmatrix
\tau(b)& U(b)\\
0                                                                 & \gamma_0(b)\\
\endsmallmatrix\right);
$$
where the intertwining matrix $U$ satisfies:
$$
U(a)=
\left(\smallmatrix
       \gp{1}\\
       \vdots\\
       \gp{1}\\
       0\\
\endsmallmatrix\right),
\quad
\quad\quad\quad
U(b)=
\left(\smallmatrix
       \gp{\alpha_1}\\
       \vdots\\
       \gp{\alpha_p}\\
       \gp{1}\\
\endsmallmatrix\right);
$$
and  $\alpha_i=\frac{\varepsilon^{p-i}-1}{\varepsilon-1}$ \quad
($i=1,2,\ldots,p$).

\proclaim {Lemma 5} $\Gamma_0$ is a faithful indecomposable
$K$-representation of $G=\gp{a,b}$.
\endproclaim

\demo{Proof} Using $\widetilde{\varepsilon}\cdot
\gp{\alpha}=\gp{\varepsilon\alpha}$ and
$1+\varepsilon+\cdots+\varepsilon^{p-1}=0$ it is easy to see that
$\Gamma_0$  is a $K$-representation. Since $\Bbb Z\subset \Bbb
Z_{(p)}\subset\Bbb Z_{p}$, it is enough to prove the lemma for
$K=\Bbb Z_p$. For this it is sufficient to prove the locality of
the centralizer of $\Gamma_0$:
$$
E(\Gamma_0)=\{\quad X \in M(p^2, K) \mid
X\cdot\Gamma_0(g)=\Gamma_0(g)\cdot X,\quad g\in G\quad \}.
$$
Let $\delta, \delta'$ be representations from  (10) and let $V$ be
a $K$-matrix such that
$$
\delta(g)\cdot V=V\cdot \delta'(g),\quad\quad (g\in G).
$$
Then
$$
V=\cases 0& \quad \text{if}\quad \delta\not=\delta';\\
 \widetilde{x},\quad \text{where}\quad x\in K[\varepsilon]& \quad
\text{if}\quad \delta=\delta'\not=\gamma_0.\\
\endcases
$$
It  follows that $X\in E(\Gamma_0)$ has the following form
$$
X=
 \left(\smallmatrix
 \widetilde{x}_1& 0&\cdots&\cdots& 0&\gp{y_1}\cr
 & \widetilde{x}_2& 0& \cdots&  0&\gp{y_2}\cr
 &&\ddots&\ddots& \vdots          & \vdots\cr
&&&\widetilde{x}_{p}&0&\gp{y_{p}}\cr
&&&&\widetilde{x}_{p+1}&\gp{y_{p+1}}\cr 0&&&&&x_0\cr
\endsmallmatrix\right),
$$
 where $x_i=x_0+(\varepsilon-1)y_i$, $x_0\in K$ and $y_i\in
K[\varepsilon]$ \quad $(i=1,2,\ldots,p+1)$. From the structure of
the matrix $X$ and the condition $K=\Bbb Z_p$ we get that $X$ is
an invertible  matrix if and only if $x_0$  is a unit in $K$.
Since $K$ is a local ring, then $E(\Gamma_0)$ is also local. The
lemma is proved.
\enddemo

\bigskip

Let $M_0=K^{p^2}$ be  the $K$-module  of the $K$-representation
$\Gamma_0$ of $G$ consisting of $p^2$-dimensional columns over
$K$. It is convenient to further condense the elements of $M_0$,
considering that the initial $p+1$ coordinates  of our vector are
$\gp{x\vert x\in K[\varepsilon]}$ (i.e. belongs to $K^{p-1}$), and
the final coordinate belong to $K$. We will do the same with
elements of $FM_0$ (the space of $p^2$-dimensional columns over
$F$).

\bigskip

\proclaim {Lemma 6} Let $\alpha=(\varepsilon-1)\m1$ and let $X, Y$
be  the following elements from $FM_0$:
$$ X=\pmatrix
       \gp{0}\\
       \vdots\\
       \gp{0}\\
       \gp{\alpha}\\
       0\\
       \endpmatrix;
\quad \quad
Y= \pmatrix
   \gp{\alpha}\\
   \vdots\\
   \gp{\alpha}\\
   \gp{0}\\
   0\\
 \endpmatrix.\tag11
$$
There  exists a $1$-cocycle $f:G=\gp{a,b}\cong C_p\times C_p\to
\widehat{M_0}= FM_0^{+}/M_0^{+}$ such that
$$
f(a)=X+M_0 \quad \text{ and }\quad f(b)=Y+M_0.
$$
Moreover,  this cocycle $f$ is not cohomologous to the zero
cocycle at each nontrivial element of $G$.
\endproclaim
\demo{Proof} Note that $\alpha=(\varepsilon-1)\m1\in \frak
F(\varepsilon)$ does not belong to $K[\varepsilon]$, but
$p\alpha\in K[\varepsilon]$. It is easy to see that the initial
$p+1$ diagonal quadratic blocks of the matrix
$$
(\Gamma_0^{p-1}+\Gamma_0^{p-2}+\cdots+\Gamma_0+E_{p^2})(g),
\quad\quad(g\in G)
$$
are either zero or have the form $p\cdot \widetilde{1}$,  and that
the final $1$-dimensional block is equal to $p$. It follows that
$$
\aligned
(\Gamma_0^{p-1}(a)+\Gamma_0^{p-2}(a)+\cdots+\Gamma_0(a)
+E_{p^2})X&\in M_0,\\
(\Gamma_0^{p-1}(b)+\Gamma_0^{p-2}(b)+\cdots+\Gamma_0(b)
+E_{p^2})Y&\in M_0,\\
(\Gamma_0(a)-E_{p^2})Y-(\Gamma_0(b)-E_{p^2})X& \in M_0.\\
\endaligned\tag12
$$
The third condition follows  from the equation
$(\widetilde{\varepsilon}- \widetilde{1})\gp{\alpha}=\gp{1}\in
K^{p-1}$.

Define a  function $f:G=\gp{a,b}\cong C_p\times C_p\to
\widehat{M_0}$ by
$$
\aligned
f(1)&=M_0; \\ f(a^i)&=(a^{i-1}+\cdots+a+1)X+M_0;\\
f(b^j)&=(b^{j-1}+\cdots+b+1)Y+M_0; \\
f(a^ib^j)&=a^if(b^{j})+f(a^i),
\endaligned
\tag13
$$
where $i,j=1,\ldots,p-1$.

According to (11)--(13) we get that $f$ is a cocycle of  $G$ with
value in   $\widehat{M_0}$. To prove the rest of the statement it
is sufficient to consider   only the generating elements
$\{a,a^ib\mid i=0,\ldots,p-1\}$ of all different nontrivial cyclic
subgroups of $G$.

If $x\in FM_0$, then by $x_{(s)}$ we define the $s^{th}$ condensed
coordinate of the vector $x$ ($s=1,\ldots,p+1$). We will do the
same  with elements of $\widehat{M_0}$. Then by (13) we obtain
that
$$
f_{(s)}(a^sb)=\gp{\varepsilon^s \alpha}+ K^{p-1}\quad \quad
(s=1,\ldots,p),\quad\quad f_{(p+1)}(a)=\gp{\alpha}+ K^{p-1}.\tag14
$$
It is easy to see that
$$
\Gamma_0(a^s)=
 \left(\smallmatrix \widetilde{\varepsilon}^s& 0&\cdots &\cdots&0& \gp{\beta_s}\cr
 & \widetilde{\varepsilon}^s& 0&\cdots&0& \gp{\beta_s}\cr
 &&\ddots&\ddots&\vdots&\vdots \cr
&&&\widetilde{\varepsilon}^s&0&\gp{\beta_s}\cr
0&&&&\widetilde{1}&0\cr &&&&&1\cr
\endsmallmatrix\right),
$$
where $\beta_s=\frac{\varepsilon^s-1}{\varepsilon-1}$
($s=1,2,\ldots,p$). Since $\varepsilon^s\alpha_s+\beta_s=0$
($s=1,2,\ldots,p$) (see the notation before Lemma 5) we have that
$(p-1)$ rows of the matrix $\Gamma_0(a^sb)-E_{p^2}$ corresponding
to the $s^{th}$ diagonal block will be zero. Besides, the final
$p$ rows of the matrix $\Gamma_0(a^sb)-E_{p^2}$ are also zero.
Thus, for any vector $z\in FM_0$ the $s^{th}$ condensed coordinate
of the vector $(\Gamma_0(a^sb)-E_{p^2})z$ $(s=1,2,\ldots,p)$ will
be equal to zero.  The $(p+1)^{th}$ coordinate in
$(\Gamma_0(a)-E_{p^2})z$ will also be zero.

Hence, according to (14) and the condition
$\alpha=(\varepsilon-1)\m1\not\in K[\varepsilon]$ it follows that
$$
(\Gamma_0(a^sb)-E_{p^2})z+f(a^sb)\not= M_0 \quad\quad \text{and}
\quad\quad (\Gamma_0(a)-E_{p^2})z+f(a)\not= M_0
$$
for any $z\in FM_0$ $(s=1,2,\ldots,p)$. The lemma is proved.

\enddemo

\proclaim {Corollary 1} The group $\cry{G;M_0;f}$ is  torsionfree.
\endproclaim

Let us define a $K$-representation of the group $G=\gp{a,b}$ as follows. Set
$$
\Delta_n=\left(\smallmatrix E_n\otimes {\gamma_3}& 0&
u_{11}& u_{12}\\
& E_n\otimes \gamma_2&
u_{21}& u_{22}\\
&&E_n\otimes \gamma_1& 0\\
0&&& E_n\otimes \gamma_0\\
\endsmallmatrix\right),
$$
where
$$
\alignat2
u_{11}(a)&=u_{21}(a)=-u_{21}(b)=E_n\otimes
\widetilde{1},\quad \quad &u_{11}(b)&=u_{22}(b)=0,\\
u_{12}(a)&=u_{12}(b)=J_n\otimes \gp{1},  \quad \quad&
u_{22}(a)&=E_n\otimes \gp{1},\\
\endalignat
$$
and $J_n$ is the upper triangular  Jordan block  of degree $n$.

\proclaim {Lemma 7}  The $K$-representation of the group $G=\gp{a,b}$
defined by
$$
a\to \Delta_n(a), \quad\quad\quad\quad b\to \Delta_n(b)
$$
is indecomposable.
\endproclaim
\demo{Proof} See \cite{1,5}.
\enddemo
\bigskip

Using the representation $\Gamma_0$ we define  the following
$K$-representation of  $G$:
$$
\Gamma_n=\left(\smallmatrix
\Gamma_0 & V_n \\
0 & \Delta_n
\endsmallmatrix\right),
$$
where $V_n$  is the matrix, which elements  are intertwining
functions of the composition factors in $\Gamma_0$ with the
composition factors in $\Delta_n$. All these intertwining
functions are zero except the function $v$ which intertwines
$\gamma_3$ in $\Gamma_0$  with  the first representation
$\gamma_0$ in $E_n\otimes \gamma_0$ and $v(a)=v(b)=\gp{1}$. Thus
$$ V_n(a)=V_n(b)=
\left(\smallmatrix
0\ldots 0\quad &0&0&\ldots&0\\
\ldots \ldots \quad&\ldots&\ldots&\ldots&\ldots\\
0\ldots 0\quad&0&0&\ldots&0\\
&&&&\\
0\ldots 0\quad&\gp{1}& \gp{0}& \cdots&\gp{0}\\
0\ldots 0\quad&\gp{0}& \gp{0}& \cdots&\gp{0}\\
0\ldots 0\quad&\gp{0}& \gp{0}& \cdots&\gp{0}\\
0\ldots 0\quad&{0}& {0}& \cdots& {0}\\
\endsmallmatrix\right).
$$

\proclaim {Lemma 8}  $\Gamma_n$ is an  indecomposable
$K$-representation of $G=\gp{a,b}$.
\endproclaim
\demo{Proof} Of course $\Gamma_n$ is a $K$-representation which
is equivalent to the following $K$-representation of the group
$G$:
$$
\Gamma_n'= \pmatrix \rho_{p-1}+\cdots +\rho_{2}&V_n'\\
0&\Delta_{n+1}'\\
\endpmatrix,\tag15
$$
where $\Delta_{n+1}'$  differs from $\Delta_{n+1}$ only by the
intertwining matrix  $U'=(u'_{ij})$ (the notation for
$\Delta_n$, $u_{ij}$ was introduced after  Corollary 1, and that for
$\Gamma_0$, before Lemma 5):
$$
\aligned
u'_{12}(a)&=u'_{12}(b)=J_{n+1}'\otimes\gp{1},\\
u'_{22}(a)&=E_{n+1}\otimes\gp{1},\\
u'_{11}(b)&=u'_{22}(b)=0,\\
u'_{11}(a)&=u'_{21}(a)=-u'_{21}(b)= \pmatrix \widetilde{0}&0\\
0& E_n\otimes \widetilde{1}\\
\endpmatrix.
\endaligned\tag16
$$
Moreover, in the representation $\Delta_{n+1}'$ there is a
non-zero intertwining  between  the first $\gamma_1$ and the first
$\gamma_0$: $u(a)=0$, $u(b)=\gp{1}$. Note that we obtained
$\Gamma'_n$ from $\Gamma_n$ by a permutation of the indecomposable
components, where intertwining functions of $\Gamma'_n$ were
obtained from the corresponding ones of $\Gamma_n$. If $\Gamma'_n$
is decomposable, then either the representations
$\rho_{p-1},\ldots,\rho_2$ or their sum cannot  be components in
$\Gamma'_n$. Each of these representations has nonzero
intertwining with $\gamma_0$, which cannot be changed without
changing the zero intertwining  for $\rho_{p-1},\ldots,\rho_2$.
Thus for the decomposability  of the representation  $\Gamma'_n$
the decomposability of the representation $\Delta'_{n+1}$  is
necessary.

The $K$-representation $\Delta'_{n+1}$ of $G$ is indecomposable.
Indeed, the additive group of the intertwining functions for any
pairs of different irreducible $K$-representations  (10) of the
group $G$ is isomorphic to the additive group of the field $K_p=K/
pK$. Any equivalence transformation (over $K$) acting on
$\Delta'_{n+1}$ will change the intertwining functions  of the
different pairs of the irreducible components of  $\Delta'_{n+1}$.
If we change the intertwining functions   by  elements of the
field $K_p$ then the effect of  the functions  is replaced by the
effect of elements of the field  $K_p$. As a consequence the
$K$-representation $\Delta'_{n+1}$ can be parametrized by the
following matrix  over $K_p$:
$$
C= \pmatrix E_n'& J_{n+1}'\\
E_n'& E_{n+1}\\
\endpmatrix,
$$
where
$E_n'= \left(\smallmatrix 0&0\\
0& E_{n}\\
\endsmallmatrix\right)$
(the notation for $\Delta_n$ was introduced before Lemma 7 and in
(15)-(16)).

 The representation $\Delta_{n+1}'$ is decomposable over $K$ if and
only if  there exist  matrices $S_i\in GL(n+1, K_p)$
($i=1,\ldots,4$) such that
$$
\pmatrix S_1& 0\\
0& S_2\\
\endpmatrix\m1\cdot
C\cdot
\pmatrix S_3& 0\\
0& S_4\\
\endpmatrix=
\pmatrix E_n'& X\\
E_n'& E_{n+1}\\
\endpmatrix,\tag17
$$
where
$$
X=
\pmatrix X_1& 0\\
0& X_2\\
\endpmatrix\tag18
$$
is decomposable over $K_p$,  and $X_1,X_2$ are square matrices.

Now, suppose that $\Delta_{n+1}'$ is decomposable and satisfies
(17)-(18). It follows that
$$
S_1=
\pmatrix t_1\m1& 0\\
*& S\\
\endpmatrix, \qquad S_2=S_4=
\pmatrix t_2& 0\\
*& S\\
\endpmatrix,
$$
where  $t_1\cdot t_2\not=0,\quad S\in GL(n, K_p)$ and
$$
X=S_1\m1\cdot J_{n+1}'\cdot S_2=T\m1\cdot Y\cdot T,\tag19
$$
where
$$
T=
\pmatrix 1& 0\\
0& S\\
\endpmatrix,
\quad\quad\quad\quad\quad Y=
\pmatrix t_1\cdot t_2& y_{12}\\
y_{21}& Y_n\\
\endpmatrix.\tag20
$$
Here  $Y$ has the following description: $t_1\cdot t_2\not=0$,
$y_{12}=(t_1,0,\ldots,0)$, $y_{21}$ is  an $n$-dimensional column,
$Y_n$ is a matrix which is obtained from the matrix $J_n'$ by
changing the first column by an $n$-dimensional column over $K_p$.
Since $y_{12}S\not=0$, then the case $X_1=t_1t_2$, $X_2=S\m1Y_nS$
is impossible, where $X_1,X_2$ are defined in (18). Thus
$$
S\m1\cdot Y_n\cdot S=\pmatrix *& 0\\
0& X_2\\
\endpmatrix.\tag21
$$
Let $\overline{K_p}$ be the algebraic closure  of the field
${K_p}$. The equivalence transformation with $T$ over
$\overline{K_p}$ given in (19) can be used to continue the
decomposition of the matrix $X$ such that $X_2$ splits into Jordan
blocks over $\overline{K_p}$ (see (17)-(21)).

Of course, we can arrange that $X_2$ is $J_s(\alpha)$, the Jordan block
with $\alpha$ in the main diagonal. The matrix $Y$ we can be considered
as a linear operator on the space ${\overline{K_p}}^{n+1}$ of
${n+1}$-dimension column vectors. Thus from (18) it follows that
$X$ is the matrix of the operator $Y$ in that basis of the space
$\overline{K_p}^{n+1}$, which consist of  the columns of the
matrix $T$. The Jordan block $X_2=J_s(\alpha)$ corresponds to the
eigenvector $e\in \overline{K_p}^{n+1}$ of the operator $Y:$\quad
$Ye=\alpha e$. Since $X_2$  does not include the first column of
$X$, the vector $e$ is a column of the matrix $T$, which is
different from the first column, i.e. the first component of $e$
is equal to zero. Using the description of $Y$ in (19), it is easy
to show that the equality $Ye=\alpha e$ ($\alpha\in
\overline{K_p}$) is impossible for the vector
$e=(0,\gamma_1,\ldots,\gamma_n)^T\not=0$. This  contradicts the
decomposability of the $K$-representation $\Delta_{n+1}'$ of $G$.
The lemma is proved.
\enddemo

\demo{Proof of Theorem 2} We can suppose that $M\subset
F^{d_n}$. Each $K$-basis in $M$ is also an $F$-basis in $F^{d_n}$
and an $(F^+/K^+)$-basis in $\widehat{M}={F^{d_n}}^+/M^+$, where
${d_n}=\deg(\Gamma_n)=(3p-2)n+p^2$. Thus the dimension of the
group $\cry{G;M_n;T_n}$ is equal to  ${d_n}$, the  degree of the
representation $\Gamma_n$, and is an unbounded function of $n$. The
theorem is proved.
\enddemo

\head
5. Proof of Theorem 3
\endhead

In  this chapter we suppose that  $K=\Bbb Z$ and we will consider
only classical crystallographic groups. Let
$A_4=\gp{\,a,b\mid a^2=b^3=(ab)^3=1\,}$ be the alternating
group of degree  $4$. Using \cite{8}, we consider the following
$\Bbb Z$-representations of $A_4$:
$$
\spreadlines{2\jot}
\alignat{6}
&\Delta_1:
&&\quad a\to1,
&&\quad b\to1;
&&\qquad \Delta_2:
&& \quad a\to\left(\smallmatrix 1&0\\ 0&1\endsmallmatrix\right),
&&\quad b\to\left(\smallmatrix 0&-1\cr1&-1 \endsmallmatrix\right);
\\
&\Delta_3:&&\quad a\to\left(\smallmatrix
                      0&-1&1\\
                      0&-1&0\\
                      1&-1&0 \endsmallmatrix\right),
&&\quad b\to \left(\smallmatrix
                      0&0&1\\
                      1&0&0\\
                      0&1&0 \endsmallmatrix\right) ;
&&\qquad\Delta_4:
&&\quad a\to \left(\smallmatrix
                      \!-1&\!-1&\!-1\cr
                      \!\phantom{-}0&\phantom{-}0&\phantom{-}1\cr
                     \!\phantom{-}0&\phantom{-}1&\phantom{-}0
\endsmallmatrix\right), &&\quad b\to \left(\smallmatrix
                     0&-1&1\cr
                     1&\phantom{-}0&1\cr
                     0&\phantom{-}1&0
\endsmallmatrix\right);
\\
&\Delta:
&&\quad a\to \Delta(a),
&&\quad b\to \Delta(b);
&&\qquad\Gamma_n:
&&\quad  a\to \Gamma_n(a),
&&\quad b\to \Gamma_n(b);
\endalignat
$$
where
$$
\Delta= \left(\smallmatrix \Delta_3&0&X_1&X_3\cr
                     &\Delta_3&X_2&0\cr
                      &&\Delta_2&0 \cr
                      0&&&\Delta_4\cr
                      \endsmallmatrix\right),\quad\quad\quad\quad\quad
                      \quad \quad
                      \Gamma_n= \left(\smallmatrix E_n\otimes \Delta_1&U\cr
                    0&E_n\otimes \Delta\cr
                      \endsmallmatrix\right),
$$
$$
X_1(a)=\left(\smallmatrix \!\!\phantom{-}1&0\cr
                            \!\!\phantom{-}0&1\cr
                           \!\!-1&1\cr
              \endsmallmatrix\right),\quad\quad\quad
X_2(a)=\left(\smallmatrix \!\!\phantom{-}0&1\cr
                            \!\!-1&1\cr
                           \!\!-1&0\cr
              \endsmallmatrix\right),\quad\quad\quad
X_3(a)=\left(\smallmatrix 0&\phantom{-}0&1\cr
                            0&\phantom{-}0&0\cr
                           0&-1&0\cr
              \endsmallmatrix\right),
$$
$$
X_i(b)=0,\quad\quad i=1,2,3,
$$
$$
U(a)=E_n\otimes \alpha+J_n(0)\otimes \beta,\quad \quad  U(b)=0,
$$
$$
\alpha=(0,0,0,0,2,0,1,-1,0,0,0), \quad\quad
\beta=(0,-2,0,0,0,0,0,1,-1,-1,0),
$$
$J_n(\nu)$ is the Jordan block  of degree $n$ with $\nu$ in the
main diagonal.

As was proved in \cite{8}, the representations
$\Delta_1,\Delta_2,\Delta_3$ and $\Delta_4$ are irreducible and
$\Delta$ and $\Gamma_n$ are indecomposable $\Bbb
Z$-representations.

Let $M_n$ be  a $\Bbb Z$-module affording the representation
$\Gamma_n$ of \, $A_4$  consisting of
$d_n$ dimensional columns over $\Bbb Z$, where
$$
\deg(\Gamma_n)=d_n=12n.
$$
It is easy to check  that
$f_n:A_4\to \widehat{M_n}$ defined by
$$
\textstyle
f_n(a)=(\underbrace{0,\ldots,0}_{n+3},\frac{1}{2},
\frac{1}{2}, 0,\ldots,0)^{T}+M_n; \quad\quad f_n(b)=(\frac{1}{3},
0,\ldots,0)^{T}+M_n
$$
is a $1$-cocycle, which is special. Therefore, we obtain
\proclaim
{Corollary 2} The classical crystallographic group
$\cry{A_4;M_n;f_n}$ is torsionfree.
\endproclaim

\enddocument